\documentclass[oneside,12pt]{amsart}
\usepackage[T1]{fontenc} 
\usepackage[a4paper]{geometry} 
\usepackage[utf8]{inputenc}
\usepackage{graphicx}

\usepackage[dvipsnames]{xcolor} 
\usepackage[pdftex,  colorlinks=true]{hyperref} 
    \hypersetup{urlcolor=RoyalBlue, linkcolor=RoyalBlue,  citecolor=black}
\usepackage{setspace} 

\usepackage{amssymb,amsthm}
\usepackage{enumitem}
\usepackage{tikz,tikz-cd}
\usepackage{url}
\usepackage[caption=false]{subfig}
\usepackage{float}
\usetikzlibrary{matrix}
\newcommand{\Ker}{\mathrm{Ker}}

\newcommand{\im}{\mathrm{Im}}

\newcommand{\PGL}{\mathrm{PGL}}
\newcommand{\PSL}{\mathrm{PSL}}

\newcommand{\ZZ}{\mathbb{Z}}
\newcommand{\RR}{\mathbb{R}}
\newcommand{\RH}{\RR\textbf{H}}
\newcommand{\QQ}{\mathbb{Q}}

\newtheorem{thm}{Theorem}[section]{\bfseries}{\itshape}
\newtheorem{theorem}[thm]{Theorem}{\bfseries}{\itshape}
\newtheorem{defn}[thm]{Definition}{\bfseries}{}
{\bfseries}{}
\newtheorem{cor}[thm]{Corollary}

\title[Cohomology of NEC Groups]{Cohomology of Fuchsian Groups and Non-Euclidean Crystallographic Groups}
  
\usepackage[foot]{amsaddr}
\author{Sam Hughes}
\address{Mathematical Institute, Andrew Wiles Building, Observatory Quarter, University of Oxford, Oxford OX2 6GG, UK}
\email{sam.hughes@maths.ox.ac.uk}
\date{\today}

\begin{document}

\maketitle

\begin{abstract}
For each geometrically finite $2$-dimensional non-Euclidean crystallographic group (NEC group), we compute the cohomology groups.  In the case where the group is a Fuchsian group, we also determine the ring structure of the cohomology.
\end{abstract}

\section{Introduction}
Let $\Gamma$ be a geometrically finite non-Euclidean crystallographic group (NEC group),  i.e. a discrete subgroup of $\PGL_2(\RR)$ with a finite sided fundamental domain for the action of $\Gamma$ on the hyperbolic plane $\RH^2$.  Throughout we let $\Lambda(\Gamma)$ denote the \emph{limit set} of $\Gamma$.  In this paper, we will calculate the cohomology of $\Gamma$.  In the case where $\Gamma$ is a Fuchsian group, i.e. $\Gamma$ is contained in $\PSL_2(\RR)$, we will also calculate the cohomology ring.  Our proof will involve finding a suitable fundamental domain for the action of the group on $\RH^2\cup \Lambda(\Gamma)$ and then applying a Cartan-Leray type spectral sequence.

Since $\RH^2\cup\Lambda(\Gamma)$ is contractible, the sequence converges to the cohomology of $\Gamma$.  Using knowledge of the abelianization of $\Gamma$, it is easy to compute with the spectral sequence.  We will now set the convention that an omission of coefficients in the (co)homology functors should be read as having coefficients in the trivial module $\ZZ$.

\begin{defn}\label{Definiton}
Let $m_1,\dots,m_r$ be a set of positive integers each greater than $2$. For $j=1,\dots,r-1$, let $\hat{t}_j$ be the greatest common divisor of the set of products of $m_1,\dots,m_r$ taken $j$ at a time.  Then, let $t_1=\hat{t_1}$ and for $j=2,\dots,r-1$ let $t_{j}=\hat{t}_j/\hat{t}_{j-1}$.  We define $w_j$ for $j=1,\dots,r-1$ by the same process except for we perform the procedure to $2m_1,m_1,\dots,m_r$ and discard products containing $2m_1m_1$.  Finally, we define $w_r$ to be equal to $2m_1m_2\dots m_r/w_{r-1}$.
\end{defn}

\begin{thm}\label{theorem.main}
Let $\Gamma$ be an NEC group of signature
\[(g,s,\epsilon,[m_1,\dots,m_r],\{(n_{1,1},\dots,n_{1,s_1}),\dots,(n_{k,1},\dots,n_{k,s_k}),(),\dots,()\}),\]
where the number of empty cycles equals $d$.  Let $C_E$ denote the number of even $n_{i,l}$ and let $C_O$ denote the number of period cycles for which every $n_{i,l}$ is odd.
\begin{enumerate}[label=(\alph*)]
    \item If $\epsilon=+$ and $d=k=s=0$  (i.e. $\Gamma$ is a cocompact Fuchsian group) then \label{theorem.main.1}
\[H_q(\Gamma) = \left\{ \begin{array}{lcr}
\ZZ && q=0,\\
\ZZ^{2g} \oplus\left(\bigoplus_{j=1}^{r-1}\ZZ_{t_j}\right) &&  q=1,\\
\ZZ &&   q=2,\\
\bigoplus_{j=1}^r\ZZ_{m_j} &&  q=2l+1,\textrm{ where } l\geq1,\\
0 && \textrm{otherwise.}
\end{array}\right.\]  

    \item If $\epsilon=-$ and $d=k=s=0$ then \label{theorem.main.4}
\[H_q(\Gamma) = \left\{\begin{array}{lcr}
\ZZ && q=0,\\
\ZZ^{g-1} \oplus\left(\bigoplus_{j=1}^{r}\ZZ_{w_j}\right) &&  q=1,\\
\bigoplus_{j=1}^r\ZZ_{m_j} &&  q=2l+1,\textrm{ where } l\geq1,\\
0 && \textrm{otherwise.}
\end{array}\right.\]
    
    \item If $\epsilon=+$ and $d+k+s>0$ then \label{theorem.main.2}
\[
H_q(\Gamma)=\left\{\begin{array}{lr}
\ZZ & q=0,\\
\ZZ^{2g+s+k+d-1}\oplus\ZZ_2^{C_E+C_O+d}\oplus\left(\bigoplus_{j=1}^r\ZZ_{m_j} \right) & q=1,\\
\ZZ_2^{\frac{1}{2}qC_E+C_O+d} & q=2p>0,\\
\ZZ_2^{\frac{1}{2}(q-1)C_E+C_O+d}\oplus\left(\bigoplus_{i=1}^k\bigoplus_{l=1}^{s_i}\ZZ_{n_{i,l}}\right)\\ \quad \oplus\left(\bigoplus_{j=1}^r\ZZ_{m_j} \right) & q\equiv3\pmod{4},\\
\ZZ_2^{\frac{1}{2}(q+1)C_E+C_O+d}\oplus\left(\bigoplus_{j=1}^r\ZZ_{m_j} \right) & q>1 \textrm{ and }\\ &q\equiv 1\pmod 4.
\end{array}\right.\]
    
    \item If $\epsilon=-$ and $d+k+s>0$ then \label{theorem.main.3}
\[
H_q(\Gamma)=\left\{\begin{array}{lr}
\ZZ & q=0,\\
\ZZ^{g+s+k+d-1}\oplus\ZZ_2^{C_E+C_O+d}\oplus\left(\bigoplus_{j=1}^{r}\ZZ_{m_j} \right) & q=1,\\
\ZZ_2^{\frac{1}{2}qC_E+C_O+d} & q=2p>0,\\
\ZZ_2^{\frac{1}{2}(q-1)C_E+C_O+d}\oplus\left(\bigoplus_{i=1}^k\bigoplus_{l=1}^{s_i}\ZZ_{n_{i,l}}\right)\\ \quad \oplus\left(\bigoplus_{j=1}^r\ZZ_{m_j} \right) & q\equiv3\pmod{4},\\
\ZZ_2^{\frac{1}{2}(q+1)C_E+C_O+d}\oplus\left(\bigoplus_{j=1}^r\ZZ_{m_j} \right) & q>1 \textrm{ and }\\ &q\equiv 1\pmod 4.
\end{array}\right.\]
\end{enumerate}
\end{thm}

In the case where $\Gamma$ is a Fuchsian group we also compute the ring structure (Theorem~\ref{theorem.Fgrp.Ring}).

\begin{defn}
We will write $\bigoplus_{j=1}^r\ZZ_{m_j} = (\bigoplus_{j=1}^{r-1}\ZZ_{t_j})\oplus(\bigoplus_{k=1}^l\ZZ_{q_k})$, where the $\bigoplus_{k=1}^l\ZZ_{q_k}$ term is decomposed via the invariant factor decomposition of finite abelian groups.  We write $\widetilde{H}^\ast(X)$ for the reduced cohomology of $X$, that is the kernel of the map induced by the inclusion of the basepoint.  Recall that $H^\ast(\ZZ_q)=\ZZ[x]/(qx)$ where $x$ has degree $2$.  Define $R_{q}$ to be the subring of $\widetilde{H}^\ast(\ZZ_q)$ generated by $x^2$ and $x^3$.
\end{defn}

\begin{thm}\label{theorem.Fgrp.Ring}
Let $\Gamma$ be a Fuchsian group of signature $[g,s;m_1,\dots,m_r]$.
\begin{enumerate}[label=(\alph*)]
\item If $s=0$ then $\widetilde{H}^\ast(\Gamma)\cong \widetilde{H}^\ast(\Sigma_g)\oplus \left(\bigoplus_{j=1}^{r-1}\widetilde{H}^\ast(\ZZ_{t_j})\right)\oplus(\bigoplus_{k=1}^l R_{q_k})$.
\item If $s>0$ then $\widetilde{H}^\ast(\Gamma)\cong \widetilde{H}^\ast(F_{2g+s-1})\oplus \left(\bigoplus_{j=1}^r\widetilde{H}^\ast(\ZZ_{m_j})\right)$ where $F_{2g+s-1}$ is a free group of rank $2g+s-1$.
\end{enumerate}
\end{thm}

We remark that some of the results have appeared in the literature before.  The case where $\Gamma$ is a cocompact Fuchsian group, so $\epsilon=+$ and $d=k=s=0$, was considered by Majumdar \cite{Majumdar1970}, however, our computation of the ring structure is new. The case $\epsilon=+$ and $d=k=0$ is a corollary of a result of Huebschmann \cite{Huebschmann1979} and the case $\epsilon=-$ and $d=k=s=0$ was considered by Akhter and Majumdar \cite{akhter2003determination}.  Each of these previous results used different methods to the ideas here. 

Other interpretations of the cohomology of Fuchsian groups have appeared in the literature.  These have primarily dealt with lifting phenomena \cite{Patterson1975}, with Eichler cohomology \cite{Curran1970,Eichler1957}  or with $K$-theory in relation to the Baum-Connes conjecture \cite{BerkovePinedaPearson2002,HughesKandKO2020,luck2000computations}.

The paper is structured as follows.  In Section~\ref{sec.NEC} we define the signature of an NEC group.  In Section~\ref{sec.Equi} we introduce the Cartan Leray type spectral sequence for a $\Gamma$-space.  Finally, in Section~\ref{coho} we prove Theorem~\ref{theorem.main} and Theorem~\ref{theorem.Fgrp.Ring}.

\section{Non-Euclidean crystallographic groups}\label{sec.NEC}
We will first describe Wilkie and Macbeath's NEC signatures \cite{Wilkie1966,Macbeath1967}, then the associated fundamental domain in $\RH^2\cup\Lambda(\Gamma)$, and finally we will give a presentation for an NEC group in terms of its signature.  For further information on NEC groups the reader should consult \cite{SOCRS2010}.

An \emph{NEC signature} consists of a sign $\epsilon=\pm$, and several sequences of integers grouped in the following manner:
\begin{enumerate}
    \item Two integers $g,s\geq0$.
    \item An ordered set of integer \emph{periods} $[m_1,\dots,m_r]$.
    \item An ordered set of $k$ \emph{period cycles} $\{C_i:=(n_{i,1},\dots,n_{i,s_i}) : 1\leq i \leq k\}$.
    \item A further $d$ empty period cycles $(),\dots,()$.
\end{enumerate}
The sequences and sign are then combined into the NEC signature, which is written as
\[(g,s,\epsilon,[m_1,\dots,m_r],\{(n_{1,1},\dots,n_{1,s_1}),\dots,(n_{k,1},\dots,n_{k,s_k}),(),\dots,()\}).\]

We let $C_E$ denote the number of even $n_{i,l}$ and we let $C_O$ denote the number of $C_i$ for which every $n_{i,l}$ is odd.

Associated to each NEC signature is a \emph{surface symbol} describing a fundamental domain for the associated NEC group.  The surface symbol is a list of edges travelling around the polygon clockwise.  Two edges paired orientably will be indicated by the same letter and a prime.  Two edges paired non-orientably will be indicated by the same letter and an asterisk.
When $\epsilon=+$, we have the surface symbol
\[\xi_1\xi_1'\dots\xi_r\xi_r'\epsilon_1\gamma_{1,0}\dots,\gamma_{1,s_1}\epsilon_1'\epsilon_2\dots\epsilon_k\gamma_{k,0}\dots,\gamma_{k,s_k}\epsilon_k'\alpha_1\beta_1'\alpha_1'\beta_1'\dots\alpha_g\beta_g'\alpha_g'\beta_g'. \]
When $\epsilon=-$, we have the surface symbol
\[\xi_1\xi_1'\dots\xi_{r+s}\xi_{r+s}'\epsilon_1\gamma_{1,0}\dots,\gamma_{1,s_1}\epsilon_1'\epsilon_2\dots\epsilon_k\gamma_{k,0}\dots,\gamma_{k,s_k}\epsilon_k'\alpha_1\alpha_1^\ast\dots\alpha_g\alpha_g^\ast. \]

\begin{figure}[H]
\centering
    \subfloat[\label{fig.FunDom.a}]{
    \includegraphics[width=0.4\textwidth]{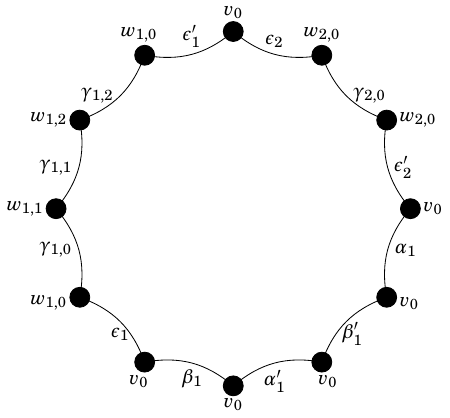}
    }
    \subfloat[\label{fig.FunDom.b}]{
    \includegraphics[width=0.4\textwidth]{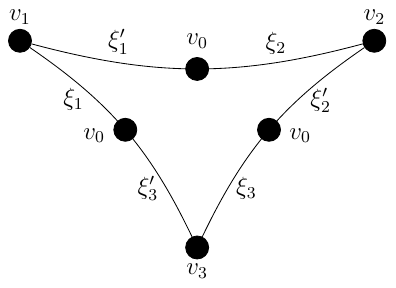}
    }\\
    \subfloat[\label{fig.FunDom.c}]{
    \includegraphics[width=0.4\textwidth]{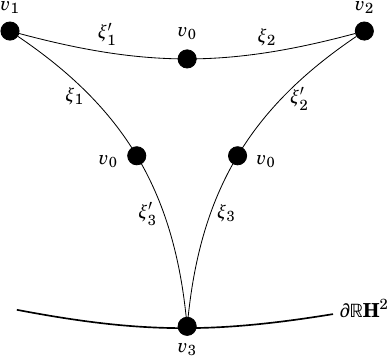}
    }
    \subfloat[\label{fig.FunDom.d}]{
    \includegraphics[width=0.4\textwidth]{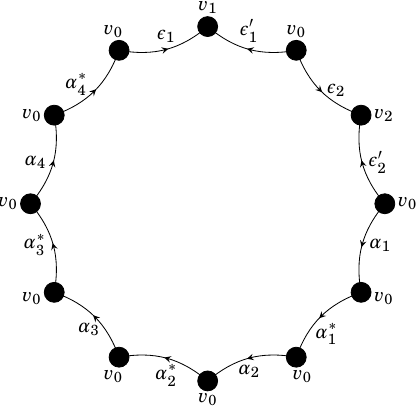}
    }
    \caption[Some fundamental domains of NEC groups.]{In (a) we have a fundamental domain for an NEC group of signature $(1,0,+,[],\{(m,n),()\})$.  The topological quotient of $\RH^2$ is homeomorphic to a torus with two open discs removed.  In the orbifold structure of the quotient we have two cone points on one of the two boundary components. 
    In (b) we have a fundamental domain for a Fuchsian triangle group of signature $(0,0,+,[p,q,r],\{\})=[0,0;p,q,r]$ for $p^{-1}+q^{-1}+r^{-1}<1$. The topological quotient is homeomorphic to a sphere.  In the orbifold structure we have three cone points.
    In (c) we have a fundamental domain of a Fuchsian NEC group of signature $(0,1,+,[m,n],\{\})=[0,1;m,n]$ for $m+n>4$.  The topological quotient is homeomorphic to a punctured sphere.  In the orbifold structure we have two cone points in the interior of the punctured sphere.
    In (d) we have a fundamental domain for a non-orientable NEC group of signature $(4,0,-,[m,n],\{\})$ for $m,n\geq2$.  The topological quotient is homeomorphic to a non-orientable surface.  In the orbifold structure we have two cone points.}
    \label{fig.FunDoms}
\end{figure}

For $j=1,\dots,r$, the period $m_j$ is attached to the vertex $v_j$ common to the edges $\xi_j$ and $\xi_j'$.  For $1\leq i\leq k$ and $1\leq l\leq s_i$ the cycle period $n_{i,l}$ is associated with the vertex $w_{i,l}$ common to the edges $\gamma_{i,l-1}$ and $\gamma_{i,l}$.  The vertices $v_j$ for $j=r+1,\dots,r+s$ lie on the boundary $\partial\RH^2$.  For $i=1,\dots,d+k$ we label the vertex common to the edges $\epsilon_i$ and $\gamma_{i,0}$ or to the edges $\gamma_{i,s_i}$ and $\epsilon_i'$ by $w_{i,0}$.  Finally, we label all other vertices $v_0$. Several examples of fundamental domains are given in Figure~\ref{fig.FunDoms}.

For an NEC group $\Gamma$ we may take the quotient $\mathcal{O}=\RH^2/\Gamma$.  The quotient comes with a natural orbifold structure and many of the geometric-toplogical features of the quotient are reflected in the signature.  Indeed, if $\epsilon=+$ then $\mathcal{O}$ is a genus $g$ surface with the disjoint union of $s$ points and $d+k$ open disks removed.  We refer to the removed points as the {\it cusps} of $\mathcal{O}$ and to the boundary of the open disks as the {\it boundary components} of $\mathcal{O}$.  There are $r$ {\it cone} or orbifold points in the interior $\mathcal{O}$.  For the $i$th boundary component, for $1\leq i\leq k$, there are $s_i$ cone or orbifold points on the boundary.  The remaining $d$ boundary components do not have any cone points.  If $\epsilon=-$ the situation is identical except we begin with a sphere with $g$ cross-caps attached.

Under the action of the associated NEC group, for $1\leq j\leq r$ the stabiliser of the vertex $v_j$ is a cyclic group of order $m_j$ acting on $\RH^2$ via rotations.  This corresponds exactly to a maximal elliptic subgroup of $\Gamma$ fixing the point $v_j$ in $\RH^2$.  If $v_j$ lies on $\partial\RH^2$, that is when $r+1\leq j\leq r+s$, then the stabiliser is isomorphic to $\ZZ$.  This corresponds to a maximal parabolic subgroup of $\Gamma$ stabilising a cusp.  

The stabiliser of the edge $\gamma_{i,l}$ for $1\leq i \leq k$ and $0\leq l\leq s_k$ or for $k+1\leq i\leq k+d$ and $l=0$ is a reflection group $\ZZ_2$.  The reflection corresponds to a non-trivial reflection in $\Gamma$ reflecting $\RH^2$ through the geodesic line containing $\gamma_{i,l}$. In the quotient these edges correspond to the edges in the boundary components. The stabiliser of the vertex $w_{i,l}$ for $1\leq i \leq k$ and $1\leq l\leq s_k$ is a dihedral group $D_{2n_{i,l}}$ of order $2n_{i,l}$, note the $w_{i,l}$ lies in the $i$th boundary component.  The stabiliser of the vertex $w_{i,l}$ for $1\leq i\leq k+d$ and $l=0$ is a reflection group $\ZZ_2$. No other points of the polygon are fixed points of the NEC group.  

Recall that the rational Euler characteristic of a group $\Gamma$ of type $VF$ is defined to be $\chi_\QQ(\Gamma)=\chi(\Gamma')/|\Gamma:\Gamma'|$ where $\Gamma'$ is a finite index subgroup of type $F$. Let $\Gamma$ be an NEC group of signature
\[(g,s,\epsilon,[m_1,\dots,m_r],\{(n_{1,1},\dots,n_{1,s_1}),\dots,(n_{k,1},\dots,n_{k,s_k}),(),\dots,()\}),\]
if $\epsilon=+$ then
\[\chi_\QQ(\Gamma)=2-2g-s-r-d-k-\frac{1}{2}\sum_{i=1}^ks_i+\sum_{j=1}^r\frac{1}{m_i}+\frac{1}{2}\sum_{i=1}^k\sum_{j=1}^{s_i}\frac{1}{n_{i,j}}\] 
and if $\epsilon=-$ then \[\chi_\QQ(\Gamma)=2-g-s-r-d-k-\frac{1}{2}\sum_{i=1}^ks_i+\sum_{j=1}^r\frac{1}{m_i}-\frac{1}{2}\sum_{i=1}^k\sum_{j=1}^{s_i}\frac{1}{n_{i,j}}.\]

If $\chi_\QQ(\Gamma)<0$ then there exists an NEC group with the corresponding signature, except when $\epsilon=-$ and $s>0$ where there is no known classification.  By the Gauss-Bonnet Theorem we see that the hyperbolic area of a fundamental domain for the NEC group is equal to $-2\pi\chi_\QQ(\Gamma)$ \cite{Singerman1974} (see also \cite[Theorem~1.1.8]{SOCRS2010}).

For the above equations, there are $17$ solutions to $\chi_\QQ(\Gamma)=0$, these exactly correspond to the $17$ Euclidean wallpaper groups \cite[Section~8]{Macbeath1967}.  We can now give a presentation for an NEC group.  Due to the large number of generators and relations, we detail this in Table~\ref{tab.NECpres}.

\begin{table}[h!]
    \centering
    \[\begin{array}{||l|l|l||}
    \hline
    \textrm{Signature element} & \textrm{Generator(s)} & \textrm{Relation(s)} \\
    \hline \hline
    \textrm{Period } m_j \textrm{ for } 1\leq j\leq r & x_j & x_j^{m_j}=1\\ \hline
    \textrm{Cycle } (n_{i,1},\dots,n_{i,s_i}) \textrm{ for } & e_i & c_{i,s_i}=e_i^{-1}c_{i,0}e_i\\
    1\leq i \leq k \textrm{ and } 0\leq l \leq s_i& c_{i,0}\dots c_{i,s_i} & c_{i,l-l}^2=c_{i,l}^2=(c_{i,l-1}c_{i,l})^2=1\\ \hline
    \textrm{Cycle } () \textrm{ for } k+1\leq i\leq k+d\ & e_i,\ c_{i,0} & c_{i,0}^2=1,\ c_{i,0}=e_i^{-1}c_{i,0}e_i\\ \hline
    s & x_r,\dots,x_{r+s} & \textrm{See } g\ \pm \\ \hline
    g\ + & a_1,b_1,\dots,a_g,b_g & \prod_{j=1}^{r+s}x_j\prod_{i=1}^{k+d}e_i\prod_{t=1}^g[a_t,b_t]=1\\ \hline
    g\ - & a_1,\dots,a_g & \prod_{j=1}^{r+s}x_j\prod_{i=1}^{k+d}e_i\prod_{t=1}^ga_t^2=1\\ \hline
    \end{array}\]
    \caption{Generators and relations for an NEC group.}
    \label{tab.NECpres}
\end{table}

If $d=k=0$ and $\epsilon=+$, then we write the signature of $\Gamma$ as $[g,s;m_1,\dots,m_r]$ and we refer to $\Gamma$ a \emph{Fuchsian group} (i.e. a discrete subgroup of $\PSL_2(\RR)$.  If $s=0$, we say that $\Gamma$ is \emph{cocompact}.

\section{A Cartan-Leray type spectral sequence}\label{sec.Equi}
For a more thorough treatment on $\Gamma$-equivariant cohomology and related spectral sequences the reader should consult for example \cite[Chapter VII]{Brown1982}.  We will just summarise the theory we need.  

Let $\Gamma$ be a discrete group, $X$ a $\Gamma$-complex and $M$ a $\Gamma$-module.  We define the {\it $\Gamma$-equivariant homology of $X$ with coefficients in $M$} to be
\[H^\Gamma_\ast(X;M):=H_\ast(\Gamma;C_\ast(X)\otimes M) \]
with diagonal $\Gamma$-action on $C_\ast(X)\otimes M$.

Let $\Omega(p)$ be a set of representatives of $\Gamma$-orbits of $p$-cells in $X$ and let $\Gamma_\sigma$ denote the stabiliser of the cell $\sigma$.  We have a $\Gamma_\sigma$-module $\ZZ_\sigma$ on which $\Gamma_\sigma$ acts on via the \emph{orientation character} $\chi_\sigma:\Gamma_\sigma\rightarrow\{\pm 1\}$.  Note that the action is trivial if $\Gamma_\sigma$ fixes $\sigma$ pointwise.  Define $M_\sigma:=\ZZ_\sigma\otimes M$, it follows $M_\sigma$ is a $\Gamma_\sigma$-module additively isomorphic to $M$ but with the $\Gamma_\sigma$-action twisted by $\chi_\sigma$.  One of the main computational tools is the following spectral sequence.

\begin{thm}\label{theorem.SpecSeq}\emph{\cite[Chapter VII (7.10)]{Brown1982}}
Let $X$ be a $\Gamma$-complex, then there is a spectral sequence \[ E^1_{pq} := \bigoplus_{\sigma\in\Omega(p)}H_q(\Gamma_\sigma;\ZZ_\sigma)\Rightarrow H^\Gamma_{p+q}(X;\ZZ). \]
\end{thm}

A description of $d^1_{p,\ast}:E_{p,*}^1\rightarrow E^1_{p-1,*}$ is given in \cite[Chapter VII.8]{Brown1982}, we will summarise it here.  Let $\sigma$ be a $p$-cell of $X$ and $\tau$ a $(p-1)$-cell.  Write $\partial_{\sigma\tau}:M_\sigma\rightarrow M_\tau$ for the $(\sigma,\tau)$-boundary component of $C_p(X)\otimes M\rightarrow C_{p-1}(X)\otimes M$.  Let $\Omega_{\sigma}=\{\tau:\partial_{\sigma\tau}\neq 0\}$ and note that this is a $\Gamma_\sigma$-invariant set of $(p-1)$-cells.  Let $\Gamma_{\sigma\tau}=\Gamma_\sigma\cap\Gamma_\tau$ and let \[t_{\sigma\tau}:H(\Gamma_\sigma;M_\sigma)\rightarrow H(\Gamma_{\sigma\tau};M_{\sigma})\] denote the transfer map arising from the fact $|\Gamma_{\sigma}:\Gamma_{\sigma\tau}|$ is finite.  Now, $\partial_{\sigma\tau}$ is a $\Gamma_{\sigma\tau}$-map and $\partial$ is a $\Gamma$-map, thus we have a map
\[u_{\sigma\tau}:H_\ast(\Gamma_{\sigma\tau};M_\sigma)\rightarrow H_\ast(\Gamma_\tau;M_\tau) \]
induced by $\Gamma_{\sigma\tau}\hookrightarrow\Gamma_\tau$ and $\partial$.  Let $\tau_0$ be a $\Gamma$-orbit representative in $X$ and choose $g\in\Gamma$ such that $g(\tau)=\tau_0$.  The action of $g$ on $C_{p-1}(X)\otimes M$ induces an isomorphism $M_\tau\rightarrow M_{\tau_0}$ which is compatible with the conjugation isomorphism $\Gamma_\tau\rightarrow\Gamma_{\tau_0}$ induced by $g$.  It follows there is an isomorphism
\[v_{\tau}:H_\ast(\Gamma_\tau;M_\tau)\rightarrow H_\ast(\Gamma_{\tau_0};M_{\tau_0}). \]
Finally, by \cite[Chapter VII (8.1)]{Brown1982} up to sign we have
\[d^1_{p,\ast}|_{H_\ast(\Gamma_\sigma;M_\sigma)}=\sum_{\tau\in\Omega(p-1)}v_\tau u_{\sigma\tau} t_{\sigma\tau}. \]

\section{Cohomology}\label{coho}
\subsection{The cocompact Fuchsian case}
We will calculate the cohomology of cocompact Fuchsian groups.  We note that the proof here is new, except for we calculate the abelianization using Smith normal form in the same way as Majumdar \cite{Majumdar1970}.

\begin{proof}[of Theorem~\ref{theorem.main}\ref{theorem.main.1}]
We will use Theorem~\ref{theorem.SpecSeq}.  In this case $X=\RH^2$ endowed with the induced cell structure from the Wilkie-Macbeath polygon.  To set up the spectral sequence we observe for each $m_j$ there is a $\Gamma$-orbit of $0$-cells, where each cell has stabiliser $\ZZ_{m_j}$.  Now, by Theorem~\ref{theorem.SpecSeq} the $E^1$-page of the spectral sequence has the form given by Figure~\ref{Fig.ss.Fuchsian1}.

\begin{figure}[ht]
\centering
\includegraphics[]{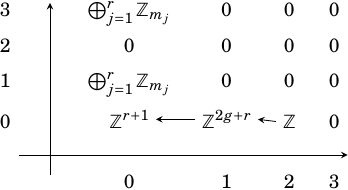}
\caption{The $E^1$-page of the spectral sequence for a Fuchsian group.}
\label{Fig.ss.Fuchsian1}
\end{figure}

The only non-trivial differentials are along the bottom row.  Slightly abusing notation we fix a basis for the chain groups by labelling the chains by the equivariant cells which afford them.  Thus, we have a sequence \[\begin{tikzcd}
0  & \arrow[l] \langle v_0,\dots,v_{r}\rangle  & \arrow[l, swap, "d^1_{1,0}"] \langle\alpha_i,\beta_i,\xi_1,\dots,\xi_r|i=1,\dots,g\rangle & \langle\gamma\rangle \arrow[l, swap, "d^1_{2,0}"] & \arrow[l] 0.\end{tikzcd}\]

We have $d^1_{1,0}(\alpha_i)=d^1_{1,0}(\beta_i)=v_0-v_0=0$ for $1\leq i\leq g$, $d^1_{1,0}(\xi_j)=v_j-v_0$ for $1\leq j\leq r$ and, $d^1_{2,0}=0$.  In particular, $\im(d^1_{1,0})\cong\ZZ^{r}$, $\Ker(d^1_{1,0})\cong\ZZ^{2g}$, $\im(d^1_{2,0})=0$ and $\Ker(d^1_{2,0})\cong\ZZ$.  From this calculation we deduce the $E^2$ page is as in Figure~\ref{Fig.ss.Fuchsian2}.

\begin{figure}[ht]
\centering
\includegraphics[]{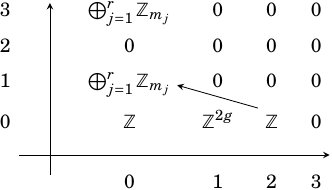}
\caption{The $E^2$-page of the spectral sequence for a Fuchsian group.}
\label{Fig.ss.Fuchsian2}
\end{figure}

The only non-trivial differential is the map drawn in Figure~\ref{Fig.ss.Fuchsian2}.  Moreover, the spectral sequence clearly collapses after the computation of this differential.  We can easily deduce what this differential is using the knowledge of $H_1(\Gamma)$.  We will compute the abelianization using the same method as Majumdar \cite{Majumdar1970}.

To compute the abelianization we write out the presentation matrix $M$ of $\Gamma$ and then compute the Smith normal form. 
\[M=\left[\begin{array}{cccccccc}
m_1    & 0      & \cdots  & \cdots  & 0             & 0      & \cdots & 0 \\
0      & \ddots &        &           & \vdots      & \vdots &       & \vdots  \\
\vdots &        & \ddots &          & \vdots      & \vdots &       & \vdots  \\
\vdots &        &        & \ddots   & 0           & 0      &       & 0  \\
0      & \cdots  & \cdots  & 0      & m_r & 0      &       & 0  \\
1      & \cdots  & \cdots  & 1      & 1   & 0      & \cdots & 0
\end{array}\right] \]
We find $H_1(\Gamma)=\ZZ^{2g}\oplus\left(\bigoplus_{j=1}^{r-1}\ZZ_{t_j}\right)$.  The constants $t_j$ (Definition~\ref{Definiton}) come from Theorem 6 in Ferrar's book `Finite Matrices' \cite{WilliamLeonardFerrar1951}.  In particular, $\prod_{j=1}^p t_j$ is equal to the greatest common divisor of the $p$-rowed minors of $M$.

It follows from the calculation of the abelianization of $\Gamma$ that the map $d^2_{2,0}$ is a surjection onto the factor $\bigoplus_{k=1}^l\ZZ_{q_k}$ from the decomposition $\bigoplus_{j=1}^r\ZZ_{m_j} = (\bigoplus_{j=1}^{r-1}\ZZ_{t_j})\oplus(\bigoplus_{k=1}^l\ZZ_{q_k})$.  The result now follows from the fact all extension problems are trivial. \qed
\end{proof}

\begin{cor}
Let $\Gamma$ be a cocompact Fuchsian group of signature $[g;m_1,\dots,m_r]$, then
\[H^q(\Gamma) = \left\{ \begin{array}{lcr}
\ZZ && q=0,\\
\ZZ^{2g} &&q=1,\\
\ZZ\oplus\left(\bigoplus_{j=1}^{r-1}\ZZ_{t_j}\right) &&  q=2,\\
\bigoplus_{j=1}^r\ZZ_{m_j} && q=2l,\textrm{ where } l\geq2,\\\
0 && \textrm{otherwise.}
\end{array}\right.\]
\end{cor}

\subsection{Non-orientable NEC groups with no cusps or boundary components}
\begin{proof}[of Theorem~\ref{theorem.main}\ref{theorem.main.4}]
Let $X=\RH^2$ and let $\Gamma$ be an NEC group with signature $(g,0,-,[m_1,\dots,m_r],\{\})$.  In this case our $E^1$-page again has the form given in Figure~\ref{Fig.ss.Fuchsian1}.  The only non-trivial differentials are along the $q=0$ row.  Keeping the same notation as before we now have a sequence 
\[\begin{tikzcd}
0  & \arrow[l] \langle v_0,\dots,v_{r}\rangle  & \arrow[l, swap, "d^1_{1,0}"] \langle\alpha_1,\dots,\alpha_g,\xi_1,\dots,\xi_r\rangle & \langle\gamma\rangle \arrow[l, swap, "d^1_{2,0}"] & \arrow[l] 0.\end{tikzcd}\]
We have $d^1_{1,0}(\alpha_i)=v_0-v_0=0$ for $1\leq i\leq g$, $d^1_{1,0}(\xi_j)=v_j-v_0$ for $1\leq j\leq r$ and, $d^1_{2,0}(f)=\sum_{i=1}^g2\alpha_i$.  In particular, $\im(d^1_{1,0})\cong\ZZ^{r}$, $\Ker(d^1_{1,0})\cong\ZZ^{2g}$, $\im(d^1_{2,0})=2\ZZ$ and $\Ker(d^1_{2,0})=0$.  It follows that $E^2_{0,0}=\ZZ$, $E^2_{1,0}=\ZZ^{g-1}\oplus\ZZ_2$ and $E^2_{2,0}=0$, the remaining entries are unchanged.  Thus, by dimension considerations the spectral sequence collapses.

The result now follows from resolving the extension problem in $H_1(\Gamma)$.  Instead we compute the abelianization of $\Gamma$ from the presentation matrix
\[M=\left[\begin{array}{cccccccc}
m_1    & 0      & \cdots  & \cdots  & 0             & 0      & \cdots & 0 \\
0      & \ddots &        &           & \vdots      & \vdots &       & \vdots  \\
\vdots &        & \ddots &          & \vdots      & \vdots &       & \vdots  \\
\vdots &        &        & \ddots   & 0           & 0      &  \cdots     & 0  \\
0      & \cdots  & \cdots  & 0      & m_r & 0      &  \cdots     & 0  \\
1      & \cdots  & \cdots  & 1      & 1   & 2      & \cdots & 2
\end{array}\right]
\sim\left[\begin{array}{cccccccc}
m_2    & 0      & \cdots  & \cdots  & 0             & 0 \\
0      & \ddots &        &           & \vdots      & \vdots  \\
\vdots &        & \ddots &          & \vdots      & \vdots  \\
\vdots &        &        & \ddots   & 0           & 0  \\
0      & \cdots  & \cdots  & 0      & m_r & 0   \\
m_1      & \cdots  & \cdots  & m_1      & m_1   & 2m_1 
\end{array}\right]=M'\]
We find $H_1(\Gamma)=\ZZ^{g-1}\oplus\left(\bigoplus_{j=1}^{r}\ZZ_{w_j}\right)$.  The constants $w_j$ (Definition~\ref{Definiton}) come from Theorem 6 in Ferrar's book `Finite Matrices' \cite{WilliamLeonardFerrar1951}.  In particular, $\prod_{j=1}^p w_j$ is equal to the greatest common divisor of the $p$-rowed minors of $M'$. \qed
\end{proof}

\begin{cor}
If $\Gamma$ is an NEC group with signature $(g,0,-,[m_1,\dots,m_r],\{\})$
then,
\[H^q(\Gamma) = \left\{\begin{array}{lcr}
\ZZ && \textrm{for }q=0,\\
\ZZ^{g-1} && \textrm{for } q=1,\\
\bigoplus_{j=1}^{r}\ZZ_{w_j} &&  \textrm{for } q=2,\\
\bigoplus_{j=1}^r\ZZ_{m_j} && \textrm{for } q=2l,\textrm{ where } l\geq2,\\
0 && \textrm{otherwise.}
\end{array}\right.\]
\end{cor}

\subsection{Orientable NEC groups with at least one cusp or boundary component}
The remaining proofs will use the homology of finite dihedral groups.  We record them here for the convenience of the reader.
\begin{theorem}\emph{\cite{handel1993}}\label{theorem.Dihedral} Let $D_{2n}$ denote a dihedral group of order $2n$.  In the case $n$ is odd we have
\[
H_q(D_{2n};\mathbb{Z}) \;=\; \left\{\begin{array}{lcl}
\ZZ && q=0,\\
\mathbb{Z}_{2} && q\equiv 1\pmod 4, \\
\mathbb{Z}_{2n} && q\equiv 3\pmod 4,\\
0 && \textrm{otherwise.} \\
\end{array}\right.
\]
\[ H_q(D_{2n};\ZZ_2)=\ZZ_2 \textrm{ for } q\geq 0. \]
In the case $n$ is even we have
\[
H_q(D_{2n};\mathbb{Z}) \;=\; \left\{\begin{array}{lcl}
\ZZ && q=0,\\
\mathbb{Z}_2^{\frac{1}{2}(q+3)} && q\equiv 1\pmod 4, \\
\mathbb{Z}_2^{\frac{1}{2}q} && q>0\textrm{ is even}, \\
\mathbb{Z}_2^{\frac{1}{2}(q+1)} \oplus \mathbb{Z}_n && q\equiv 3\pmod 4.\\
\end{array}\right.
\]
\[ H_q(D_{2n};\ZZ_2)=\ZZ_2^{q+1} \textrm{ for } q\geq 0. \]
\end{theorem}

We will now compute the cohomology of an NEC group with orientable quotient space with at least one boundary component or cusp.

\begin{proof}[of Theorem~\ref{theorem.main}\ref{theorem.main.2}]
First, assume that $k+d=0$, so $s>0$.  In this case it is easy to see that we can rearrange the presentation of $\Gamma$ so that $\Gamma\cong F_{s-1}\ast\ZZ_{m_1}\ast\dots\ast\ZZ_{m_r}$ where $F_{s-1}$ is a free group of rank $s-1$.  The result now follows from a straightforward application of the homology Mayer-Vietoris sequence.

We now treat the case with boundary, let $k,d,s\geq0$ such that $k+d>0$ and let $\epsilon=+$.  We will use Theorem~\ref{theorem.SpecSeq}; here our space $X$ is $\RH^2\cup\Lambda(\Gamma)$ endowed with the induced cell structure from the Wilkie-Macbeath polygon.  To set up the sequence, observe that the stabiliser of a marked point $v_j$ in the interior of the quotient space is a cyclic group $\ZZ_{m_j}$.  If the vertex $v_j$ lies on $\partial\RH^2$ then the stabiliser is $\ZZ$.  The stabiliser of a marked point $w_{i,l}$ on the boundary of the quotient space is a dihedral group $D_{2n_{i,l}}$, and edges along the boundary are stabilised by reflection groups isomorphic to $\ZZ_2$.  Now, the face stabilisers are trivial, vertices have a canonical orientation, and the edges being stabilised by $\ZZ_2$ are fixed pointwise.  In particular, for each cell $\sigma\in X$ the orientation character $\chi_\sigma:\Gamma_\sigma\to\{\pm1\}$ is trivial. It follows that the $E^1$-page consists of modules with trivial $\Gamma$-action and has the form given in Figure~\ref{fig.ss.NEC}.

\begin{figure}[ht]
\centering
\includegraphics[]{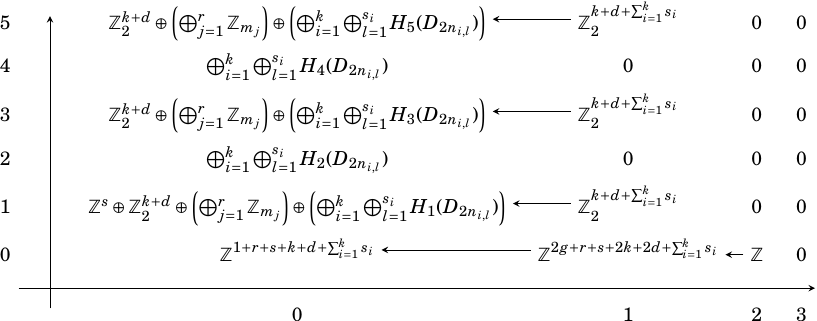}
\caption{The $E^1$-page of the spectral sequence for an orientable NEC group with cusps and boundary.}
\label{fig.ss.NEC}
\end{figure}

We will first deal with the differentials $d_{\ast,0}^1$.  By slightly abusing notation and labelling our basis elements for each chain group by the equivariant cells which afford them, we have a sequence

\[\begin{tikzcd}
0 \leftarrow \left\langle v_j, w_{i,l}\ \bigg|\ \begin{array}{c} 0\leq j\leq r+s\\ 1\leq i\leq k+d\\ 0\leq l\leq s_i\end{array} \right\rangle  & {{\Bigg\langle \alpha_t, \beta_t, \xi_j, \gamma_{i,l},\epsilon_i\ \bigg|\ \begin{array}{c} 1\leq t\leq 2g,\\ 1\leq j\leq r+s\\ 1\leq i\leq k+d,\\ 0\leq l\leq s_i\end{array} \Bigg\rangle}} \arrow[l, "{d^1_{1,0}}"'] \\
\langle f\rangle \arrow[ru, "{d^1_{2,0}}"] & 0. \arrow[l]                  
\end{tikzcd}\]

Computing the image of the differential $d_{2,0}^1$ on the $\ZZ$-basis element $f$, we obtain that up to sign \[f \mapsto \sum_{i=1}^k\sum_{l=0}^{s_i}\gamma_{i,l}.\]  So, we find $\im(d_{2,0}^1)=\ZZ$ and $E^2_{2,0}=0$.  In light of the description of the fundamental domain in Section~\ref{sec.NEC}, for $d_{1,0}^1$ we have the following
\[\begin{array}{cllcr}
    \alpha_t&\mapsto& v_0-v_0=0 & & \textrm{for } 1\leq t\leq 2g;\\
    \beta_t&\mapsto& v_0-v_0=0 & & \textrm{for } 1\leq t\leq 2g;\\
    \xi_j &\mapsto& v_j-v_0 & & \textrm{for } 1\leq j\leq r+s;\\
    \gamma_{i,l} &\mapsto& w_{i,(l+1\pmod{s_i})}-w_{i,l}\quad & & \textrm{for } 1\leq i\leq k,\textrm{ and } 0\leq l\leq s_i;\\
    \gamma_{i,0} &\mapsto& w_{i,0}-w_{i,0}=0 & & \textrm{for } k+1\leq i\leq k+d;\\
    \epsilon_i &\mapsto& w_{i,0}-v_0 & & \textrm{for } 1\leq i\leq k+d.
\end{array}\]
In particular, we have $\im(d_{1,0}^1)=\ZZ^{r+s+k+\sum_{i=1}^k s_i}$ and $\Ker(d_{1,0}^1)=\ZZ^{2g+k+d}$.  It then follows that $E^2_{1,0}=\ZZ^{2g+k+d-1}$ and $E^2_{0,0}=\ZZ$.  At this point, it is easy to see that the spectral sequence will collapse trivially once we have computed the differentials $d_{1,*}^1$.

We will begin with the differential $d_{1,q}^1$
where $q\equiv 1\pmod{4}$.  Since the edges connected to the vertices corresponding to the $\ZZ_{m_j}$ summands have trivial stabilisers, the $\ZZ_{m_j}$ summands will survive to the $E^2$-page.  In the case $q=1$, the $\ZZ$ summands also survive by the same reasoning.

We now draw our focus to the other summands.  Let each $D_{2n_{i,l}}$ be generated by a reflection $r_{i,l}$ and a rotation $t_{i,l}$ of order $n_{i,l}$.  We have that $H_1(D_{2n_{i,l}})$ is generated by $r_{i,l}^1,t_{i,l}^1$, the images of $t_{i,l}$ and $r_{i,l}$ under the abelianization map.  For $q>1$ we have classes $r_{i,l}^q,t_{i,l}^q\in H_q(D_{2n_{i,l}})$. There will be extra generators of $H_q(D_{2n_{i,l}})$ whenever an $n_{i,l}$ is even; we will suppress this from the notation.  Note that this has the effect of working in the subgroup of $H_q(D_{2n_{i,l}})$ generated by $r_{i,l}^q$ and $t_{i,l}^q$. Also, note that if $n_{i,l}$ is odd then $t^q_{i,l}=0$.  For each $q\equiv1\pmod4$ we now have a sequence (ignoring the extra classes arising from dihedral groups where $n_{i,l}$ is even and when $q>1$)
\[\begin{tikzcd}
0&\arrow[l] \left\langle w_{i,0}^q, w_{p,0}^q,r_{i,l}^q, t_{i,l}^q \ \Bigg|\ \begin{array}{c}1\leq i \leq k\\ 1\leq l \leq s_i\\ 1 \leq p\leq d \end{array} \right\rangle &\arrow[l, swap, "d_{1,q}^1"] \left\langle\gamma^q_{i,l},\gamma^q_{p,0}\Bigg| \ \begin{array}{c}1\leq i \leq k\\ 0\leq l \leq s_i\\ 1 \leq p\leq d \end{array} \right\rangle.
\end{tikzcd}\]

We will break the map $d^1_{1,q}$ into several cases depending on the adjacent edges in the fundamental domain and the cycle type of the boundary component.  First, we will consider each `end' of the $i$th boundary component with a non-empty period of cycles (i.e. $1\leq i \leq k$), the reader should keep Figure~\ref{fig.FunDom.a} in mind.  Here we have
\[(\psi_{i,0})_q:H_q(\langle r_{i,0}\rangle)\hookrightarrow H_q(D_{2n_{i,l}})\oplus H_q(\ZZ_2) \textrm{ by } \gamma_{i,0}^q\mapsto t_{i,1}^q-w_{i,0}^q \]
and
\[(\psi_{i,s_i})_q:H_q(\langle r_{i,s_i}t_{i,s_i}\rangle)\hookrightarrow H_q(\ZZ_2) \oplus H_q(D_{2n_{i,l}}) \textrm{ by } \gamma_{i,0}^q\mapsto w_{i,0}^q - t_{i,s_i}^q - r_{i,s_i}^q. \]
For the intermediary edges we have
\[(\psi_{i,l})_q:H_q(\langle r_{i,l}t_{i,l}\rangle)\hookrightarrow H_q(D_{2n_{i,l+1}}) \oplus H_q(D_{2n_{i,l}}) \textrm{ by } \gamma_{i,0}^q\mapsto  t_{i,l+1}^q - t_{i,l}^q - r_{i,l}^q. \]
In each case we are suppressing from the image a possible sum of order 2 classes (distinct from $t_{i,l}^q$ and $r_{i,l}^q$) arising from even dihedral groups.  The reason for this is that provided at least one of the $n_{i,l}$ are even, the images of the maps $\psi_{i,l}$ for $0\leq l\leq s_i$ are already linearly independent.  Of course if all of the $n_{i,l}$ for $0\leq l\leq s_i$ are odd, then the classes do not exist.

If the boundary component $i$ only contains odd cycles, then $\gamma^q_{i,s_i}=\sum_{l=0}^{s_i-1}\gamma^q_{i,l}$, so we have an order $2$ element in the kernel of $d^1_{1,q}$.  If the boundary component has an empty period of cycles, then we have exactly one edge $\gamma_{i,0}$ with vertex $w_{i,0}$ at each end.  In particular $\gamma_{i,0}^q\mapsto w_{i,0}^q-w_{i,0}^q=0$.  From this analysis we deduce that $\Ker(d^1_{1,q})=\ZZ^{C_O+d}_2$ and $\im(d^1_{1,q})\cong\ZZ_2^{k+\sum_{i=1}^k s_i-C_O}$.  It then follows from a simple calculation that $E^2_{1,q}=\ZZ_2^{C_O+d}$ and $E^1_{0,q}\cong\ZZ_2^{\frac{1}{2}(q+1)C_E+C_O+d}\oplus\left(\bigoplus_{j=1}^r\ZZ_{m_j}\right)$ for $q\equiv 1 \pmod4$, $q>1$.  When $s>0$ we have an additional $\ZZ^s$ summand in $E^2_{0,1}$.

An alternative way of considering these maps is as follows.  Let $C_{E_i}$ denote the number of even periods in the $i$th period cycle.  Observe that each period cycle contributes $\frac{1}{2}(q+1)C_{E_i}-1$ summands of $\ZZ_2$ to $E_{0,q}^2$.  The $C_O$ summands of $\ZZ_2$ contained in $\Ker(d_{1,q}^1)$ cause an additional $C_O$ summands of $\ZZ_2$ to survive to $E^2_{0,q}$. From, above we then have that \[k+\sum_{i=1}^k\left(\frac{1}{2}(q+1)C_{E_i}-1\right)+C_O=k+\frac{1}{2}(q+1)C_E-k+C_O=\frac{1}{2}(q+1)C_E+C_O.\]

We now need to compute the maps $d^1_{1,q}$ for $q\equiv 3\pmod4$.  We have essentially the same cases and proof as when $q\equiv1\pmod4$ except that $\textrm{Coker}(d^1_{1,q})$ contains a summand $\ZZ_{n_i,l}$ for each $n_{i,l}$.  

\textbf{Claim} For $q\equiv 3\pmod4$, $1\leq i\leq k$ and $1\leq l\leq s_i$ the term $E^2_{1,q}=\textrm{Coker}(d^1_{1,q})$ contains a summand $\ZZ_{n_i,l}$.

\textbf{Proof of claim:}
When $n_{i,l}$ is odd this is immediate.  Let $n:=n_{i,l}$ be even and consider $H^{q+1}(D_{2n};\ZZ)$ where $q\equiv3\pmod4$.  There is an element of order $n$ in $H^{q+1}(D_{2n};\ZZ)$ that corresponds to a power of the second Chern class of the faithful $2$-dimensional linear representation $\rho$ of $D_{2n}=\langle r,t\rangle$.  Restricting $\rho$ to the subgroup $\langle rt \rangle$ gives the regular representation of $\ZZ_2\cong\langle rt \rangle$.  Now, the total Chern class of $\ZZ_2$ is equal to $0$ in degree $4$.  It follows that the map $H^{q+1}(D_{2n};\ZZ)\rightarrow H^{q+1}(\langle rt\rangle)$ has kernel containing a $\ZZ_n$ summand. Dualizing back to homology, it follows the map $H_q(\langle rt\rangle)\rightarrow H_q(D_{2n})$ has cokernel containing a $\ZZ_n$ summand.  This yields the claim.

We conclude the description of $E^2$ as follows.  First, when $q\equiv3\pmod4$ we have $\Ker(d^1_{1,q})\cong\ZZ_2^{C_O+d}$ and $\im(d^1_{1,q})\cong\ZZ_2^{k+\sum_{i=1}^ks_i-C_O}$.  It follows $E^2_{1,q}\cong\ZZ^{C_O+d}$ and $E^2_{0,q}\cong\ZZ_2^{\frac{1}{2}(q-1)C_E+C_O+d}\oplus\left(\bigoplus_{i=1}^k\bigoplus_{l=1}^{s_i}\ZZ_{n_{i,l}}\right)\oplus\left(\bigoplus_{j=1}^k\ZZ_{m_j}\right)$.  Every other entry on the $E^2$-page is $0$ trivially.  

The theorem follows from resolving the extension problems $0\rightarrow E^2_{1,q-1}\rightarrow H_q(\Gamma)\rightarrow E^2_{0,q}\rightarrow 0$, where $q>0$ is even.  To resolve the extension problems, we will compute the homology of $\Gamma$ with $\ZZ_2$ coefficients and then compare the $\ZZ_2$-rank of $H_q(\Gamma;\ZZ_2)$ with the $\ZZ_2$-rank of $(E^2_{1,q-1}\oplus E^2_{0,q})\otimes\ZZ_2 \oplus \textrm{Tor}(E^2_{0,q-1},\ZZ_2)$.  Note that the latter is equal to $(q+1)C_E+2C_O+2d$.  If the ranks are equal, then the extension will split.

Recall that $H_n(\ZZ_2;\ZZ_2)=\ZZ_2$ for $n\geq0$.  Combining this with the $\ZZ_2$-homology groups of the Dihedral groups (Theorem~\ref{theorem.Dihedral}) and the $\Gamma$-equivariant spectral sequence (Theorem~\ref{theorem.SpecSeq}), we can set up a spectral sequence calculation.  To simplify things, note we are only interested in the maps $d^1_{1,q}$ for $q>0$.

Let $q>0$ and let $C_T$ denote the number of odd cycles, so $C_T+C_E=\sum_{i=1}^ks_i$.  We then have a sequence 
\[\begin{tikzcd}
0&\arrow[l] \ZZ_2^{(q+1)C_E+C_T+d+k} &\arrow[l, "d_{1,q}^1",swap]\ZZ_2^{C_E+C_T+d+k} & \arrow[l] 0.
\end{tikzcd}\]
By essentially using the same calculations as above we have that $\im(d_{1,q}^1)\cong\ZZ_2^{C_E+C_T+k-C_O}$.  From this we conclude that $E_{0,q}^2=\ZZ_2^{(q+1)C_E+C_O+d}$ and that $E_{1,q}^2=\ZZ_2^{C_O+d}$.  This gives a $\ZZ_2$-rank of $(q+1)C_E+2C_O+2d$.  Thus, the extension splits. \qed
\end{proof}

\begin{cor}
Let $d+k+s>0$. If $\Gamma$ is an  NEC group with signature
\[(g,s,+,[m_1,\dots,m_r],\{(n_{1,1},\dots,n_{1,s_1}),\dots,(n_{k,1},\dots,n_{k,s_k}),(),\dots,()\})\]
then,
\[
H^q(\Gamma)=\left\{\begin{array}{lr}
\ZZ & q=0,\\
\ZZ^{2g+s+k+d-1} & q=1,\\
\ZZ_2^{\frac{1}{2}qC_E+C_O+d}\oplus\left(\bigoplus_{j=1}^r\ZZ_{m_j} \right) & q\equiv2\pmod{4},\\
\ZZ_2^{\frac{1}{2}(q-1)C_E+C_O+d} & q=2p+1 \textrm{ where } p\geq 1,\\
\ZZ_2^{\frac{1}{2}qC_E+C_O+d}\oplus\left(\bigoplus_{i=1}^k\bigoplus_{l=1}^{s_i}\ZZ_{n_{i,l}}\right)\\ \quad \oplus\left(\bigoplus_{j=1}^r\ZZ_{m_j} \right) & q>0 \textrm{ and } q\equiv 0\pmod 4.
\end{array}\right.\]
\end{cor}

\subsection{Non-orientable NEC groups with at least one cusp or boundary component}
We will now compute the cohomology of an NEC group with non-orientable quotient space and at least one cusp or boundary component.  The proof is almost exactly the same as the proof of Theorem~\ref{theorem.main}\ref{theorem.main.2} so we will only provide a brief sketch and highlight the differences.

\begin{proof}[of Theorem~\ref{theorem.main}\ref{theorem.main.3} (sketch)]
First assume that $k+d>0$. The key differences between the orientable (Figure~\ref{fig.ss.NEC}) and non-orientable cases is the $E^1_{1,0}$ term and the map $d^1_{2,0}$.  The $E^1_{1,0}$ now contains a $\ZZ^g$ summand instead of a $\ZZ^{2g}$ summand.  The map $d^1_{2,0}$ now sends the generator to the sum of boundary components plus $2$ times each generator of the aforementioned $\ZZ^g$ summand.  More precisely (with the same notation as in the proof of Theorem~\ref{theorem.main}\ref{theorem.main.2}) we have,
\[ f \mapsto \sum_{i=1}^k\sum_{l=0}^{s_i}\gamma_{i,l}+2\sum_{t=1}^g\alpha_t\]
In particular, $E^2_{1,0}=\ZZ^{g+k+d-1}$.  The proof goes through identically from here. 

Now assume $g>0$ and $k+d=0$, so $s>0$.  We still have that $E^1_{1,0}$ contains a $\ZZ^g$ summand instead of a $\ZZ^{2g}$ summand.  However, with notation as before,
\[ d_{2,0}^1(f) = 2\sum_{t=1}^g\alpha_t. \]
In particular, $E^2_{1,0}=\ZZ^{g-1}\oplus\ZZ_2$.  The remainder of the proof is identical, except we now have an extension problem to determine $H_1(\Gamma)$.  We instead resolve this by computing the abelianisation from the presentation matrix
\[M=\left[\begin{array}{ccccccccccc}
m_1    & 0      & \cdots  & \cdots  & 0             & 0      & \cdots & 0 & 0& \cdots & 0 \\
0      & \ddots &        &           & \vdots      & \vdots &       & \vdots & \vdots & &\vdots  \\
\vdots &        & \ddots &          & \vdots      & \vdots &       & \vdots &\vdots & &\vdots \\
\vdots &        &        & \ddots   & 0           & 0      & \cdots      & 0 & 0 &\cdots & 0  \\
0      & \cdots  & \cdots  & 0      & m_r       &  0      & \cdots      & 0 & 0 & \cdots &0  \\
1      & \cdots  & \cdots  & 1      & 1          & 1      & \cdots & 1 & 2 & \cdots & 2
\end{array}\right]. \]
Clearly, $M$ can be reduced to an $(r+1)\times(r+g+s)$ matrix with the only non-zero entries equal to $m_1,\dots,m_r,1$ on the leading diagonal.  The result follows. 

The final case when $g=k=d=0$ and $s>0$ follows an almost identical argument to the case $k=d=0$, $s>0$ and $\epsilon=+$, so we will not recreate it here.
\qed
\end{proof}

\begin{cor}
Let $d+k+s>0$.  If $\Gamma$ is an NEC group of signature
\[(g,s,-,[m_1,\dots,m_r],\{(n_{1,1},\dots,n_{1,s_1}),\dots,(n_{k,1},\dots,n_{k,s_k}),(),\dots,()\})\]
then
\[
H^q(\Gamma)=\left\{\begin{array}{lr}
\ZZ & q=0,\\
\ZZ^{g+s+k+d-1} & q=1,\\
\ZZ_2^{C_E+C_O+d}\oplus\left(\bigoplus_{j=1}^{r}\ZZ_{m_j}\right) & q=2,\\
\ZZ_2^{\frac{1}{2}(q-1)C_E+C_O+d} & q=2p+1 \textrm{ where } p\geq1,\\
\ZZ_2^{\frac{1}{2}qC_E+C_O+d}\oplus\left(\bigoplus_{i=1}^k\bigoplus_{l=1}^{s_i}\ZZ_{n_{i,l}}\right)\\ \quad \quad \oplus\left(\bigoplus_{j=1}^r\ZZ_{m_j} \right) & q>0 \textrm{ and } q\equiv0\pmod4,\\
\ZZ_2^{\frac{1}{2}qC_E+C_O+d}\oplus\left(\bigoplus_{j=1}^r\ZZ_{m_j} \right) & q>2 \textrm{ and } q\equiv 2\pmod 4.
\end{array}\right.\]
\end{cor}

\subsection{The ring structure}
We will now deal with the computation of the ring structure.  Recall that $R_{q}$ is the subring of $\widetilde{H}^\ast(\ZZ_q)$ generated by $x^2$ and $x^3$, where $x$ is the degree $2$ generator of ${H}^\ast(\ZZ_q)$.

\begin{proof}[of Theorem~\ref{theorem.Fgrp.Ring}]
We first prove the result when $s>0$.  Let $\Gamma$ be a Fuchsian group of signature $[g,s;m_1,\dots,m_r]$ such that $s>0$.  We may rearrange the presentation of $\Gamma$ so that $\Gamma\cong F_{s-1}\ast\ZZ_{m_1}\ast\dots\ast\ZZ_{m_r}$ where $F_{s-1}$ is a free group of rank $s-1$.  The result is now an easy application of the Mayer-Vietoris sequence for cohomology.

\begin{figure}[h]
    \centering
    \includegraphics[]{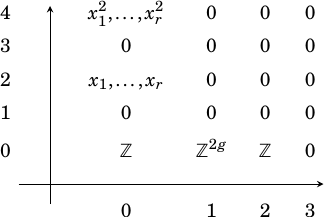}
    \caption{The $E^2$-page of the cohomological spectral sequence for a cocompact Fuchsian group.  Here the element $x_j$ is additive torsion of order $m_j$.}
    \label{fig.5}
\end{figure}

Now, let $\Gamma$ be a Fuchsian group of signature $[g,s;m_1,\dots,m_r]$ such that $s=0$.  We instead consider the equivariant cohomology spectral sequence for $\Gamma$.  Armed with our calculation for homology, it is clear that $E^2$-page has the form given in Figure~\ref{fig.5} (here $m_jx_j=0$).  Now, there is an extension problem 
\[\begin{tikzcd}
0 \arrow[r] & \ZZ \arrow[r] & H^2(\Gamma) \arrow[r] & {\langle x_1,\dots,x_r\rangle} \arrow[r] & 0
\end{tikzcd}\] which we can resolve using the computation of $H_1(\Gamma)$.  In particular, an application of the universal coefficient theorem yields that $H^2(\Gamma)\cong\ZZ\oplus\left(\bigoplus_{j=1}^{r-1}\ZZ_{t_j}\right)$.  It follows the extension problem kills a subgroup of $\langle x_1,\dots,x_r\rangle$ isomorphic to $\bigoplus_{k=1}^l\ZZ_{q_k}$.  Since the spectral sequence preserves cup products the result follows. \qed
\end{proof}

\section{Closing remarks}
We end with three remarks.  Firstly, the author was asked by Professor Gareth Jones whether the same results hold for the $17$ wallpaper groups if one takes $X$ to be the Euclidean plane.  We confirm here it does, however the cohomology computations of these groups are well known so we will not elaborate on this.  Secondly, the results in this paper are consistent with Gaboriau's result that $L^2$-Betti numbers of lattices in a Lie group are proportional to their covolume \cite{Gab02}.  As such one deduces the well known result that for an NEC group $\Gamma$ the first $L^2$-Betti number $b_1^{(2)}(\Gamma)=-\chi_\QQ(\Gamma)$ and all other $L^2$-Betti numbers vanish.  Finally, Fuchsian groups are not determined by their cohomology.  Indeed, the groups with signatures $[g,s;3,10]$ and $[g,s;5,6]$ have isomorphic cohomology rings but the groups are not isomorphic.

\section*{acknowledgements}
I would like to thank my PhD supervisor Professor Ian Leary for his guidance and support.  I would also like to thank the anonymous reviewer as their feedback greatly improved the exposition of this paper.  This work was supported by the Engineering and Physical Sciences Research Council grant number 2127970.

\end{document}